\documentclass[12pt]{article}
\usepackage{amsfonts}
\usepackage{amsmath}
\usepackage{amssymb}

\newtheorem{theorem}{Theorem}[section]

\newtheorem{conjecture}[theorem]{Conjecture}

\newtheorem{lemma}[theorem]{Lemma}

\newtheorem{proposition}[theorem]{Proposition}

\input{tcilatex}

\begin{document}

\title{A simple proof of the algebraic version of a conjecture by Vogan}
\author{Tim Bratten and Sergio Corti \\
%EndAName
Facultad de Ciencias Exactas, UNICEN. Tandil, Argentina. }
\date{}
\maketitle

\begin{abstract}
In a recent manuscript, D.Vogan conjectures that four canonical
globalizations of Harish-Chandra modules commute with certain $\mathfrak{n}-$%
cohomology groups. In this article we prove that Vogan's conjecture holds
for one of the globalizations if and only if it holds for the dual.  Using a
previously published result of one of the authors, which establishes the
conjecture for the minimal globalization, we can therefore deduce Vogan's
conjecture for the maximal globalization
\end{abstract}

\section{Introduction}

\noindent In a recent manuscript \cite[Conjecture 10.3]{V1}, D.Vogan
conjectures that four canonical globalizations of Harish-Chandra modules
commute with certain $\mathfrak{n}-$cohomology groups. In this article we
consider an algebraic version of Vogan's conjecture which entails that the
conjecture holds for one of the globalizations if and only if it holds for
the dual. We prove this result for a reductive Lie group of Harish-Chandra
class. Using the result that Vogan's conjecture is known for the minimal
globalization \cite{Br1} we can therefore conclude the conjecture is true
for the maximal globalization.

\medskip

\noindent This article is organized as follows. In the second section we
define the $\mathfrak{n}-$homology and cohomology groups and recall some
formulas that relate them. The third section treats Harish-Chandra modules
and globalizations. The fourth section introduces Vogan's conjecture and
establishes the algebraic version. In the fifth section we review the
Hochschild-Serre spectral sequence and use it to prove the main result.

\medskip

\noindent The authors would like to take this opportunity to thank David
Vogan for contacting us about his work relating to this conjecture. We would
also like to thank the referee for pointing out the possibility of a formal
proof based on the Hochschild-Serre spectral sequence. Our original proof,
which applied to the case of a connected, complex reductive group, relied on
the Beilinson-Bernstein classification and characterized the Harish-Chandra
dual of a standard module.

\section{$\mathfrak{n}-$homology and cohomology}

\noindent Suppose $G_{0}$ is a reductive group of Harish-Chandra class with
Lie algebra $\mathfrak{g}_{0}$ and complexified Lie algebra $\mathfrak{g}$.
By definition, \emph{a Borel subalgebra }of $\mathfrak{g}$ is a maximal
solvable subalgebra and \emph{a parabolic subalgebra }of $\mathfrak{g}$ is a
subalgebra that contains a Borel subalgebra. If $\mathfrak{p}\subseteq 
\mathfrak{g}$ is a parabolic subalgebra then \emph{the nilradical} $%
\mathfrak{n}$ of $\mathfrak{p}$ is the largest solvable ideal in $\left[ 
\mathfrak{p},\mathfrak{p}\right] $. A\emph{\ Levi factor }is a complementary
subalgebra to $\mathfrak{n}$ in $\mathfrak{p}$. One knows that Levi factors
exist and that they are exactly the subalgebras which are maximal with
respect to being reductive in $\mathfrak{p}$.

\medskip

\noindent Fix a parabolic subalgebra $\mathfrak{p}$ with nilradical $%
\mathfrak{n}$ and Levi factor $\mathfrak{l}$. Let $U(\mathfrak{n})$ denote
the enveloping algebra of $\mathfrak{n}$ and let $\mathbb{C}$ be the
irreducible trivial module. If $M$ is a $\mathfrak{g}-$module then \emph{the
zero} $\mathfrak{n}-$\emph{homology of} $M$ is the $\mathfrak{l}-$module 
\begin{equation*}
H_{\text{0}}(\mathfrak{n},M)=\mathbb{C}\otimes _{U(\mathfrak{n})}M.
\end{equation*}%
This definition determines a right exact functor from the category of $%
\mathfrak{g}-$modules to the category of $\mathfrak{l}-$modules. \emph{The} $%
\mathfrak{n}-$\emph{homology groups of} $M$ are the $\mathfrak{l}-$modules
obtained as the corresponding derived functors.

\smallskip

\noindent \emph{The right standard resolution of }$\mathbb{C}$ is the
complex of free right $U(\mathfrak{n})-$modules given by 
\begin{equation*}
\cdots \rightarrow \Lambda ^{p+1}\mathfrak{n}\otimes U(\mathfrak{n}%
)\rightarrow \Lambda ^{p}\mathfrak{n}\otimes U(\mathfrak{n})\rightarrow
\cdots \rightarrow \mathfrak{n}\otimes U(\mathfrak{n})\rightarrow U(%
\mathfrak{n})\rightarrow 0\text{.}
\end{equation*}%
Applying the functor 
\begin{equation*}
-\otimes _{U(\mathfrak{n})}M\text{ }
\end{equation*}%
to the standard resolution we obtain a complex 
\begin{equation*}
\cdots \rightarrow \Lambda ^{p+1}\mathfrak{n}\otimes M\rightarrow \Lambda
^{p}\mathfrak{n}\otimes M\rightarrow \cdots \rightarrow \mathfrak{n}\otimes
M\rightarrow M\rightarrow 0
\end{equation*}%
of left $\mathfrak{l}-$modules called \emph{the standard} $\mathfrak{n}-$%
\emph{homology complex}. Here $\mathfrak{l}$ acts via the tensor product of
the adjoint action on $\Lambda ^{p}\mathfrak{n}$ with the given action on $M$%
. Since $U(\mathfrak{g})$ is free as $U(\mathfrak{n})-$module, a routine
homological argument identifies the pth$-$homology of the standard complex
with the pth $\mathfrak{n}-$homology group 
\begin{equation*}
H_{\text{p}}(\mathfrak{n},M).
\end{equation*}%
One knows that the induced $\mathfrak{l}-$action on the homology groups of
the standard complex is the correct one.

\smallskip

\noindent The \emph{zero} $\mathfrak{n}-$\emph{cohomology }of a $\mathfrak{g}%
-$module $M$ is the $\mathfrak{l}-$module 
\begin{equation*}
H^{0}(\mathfrak{n},M)=\text{Hom}_{U(\mathfrak{n})}(\mathbb{C},M).
\end{equation*}%
This determines a left exact functor from the category of $\mathfrak{g}-$%
modules to the category of $\mathfrak{l}-$modules. By definition, \emph{the} 
$\mathfrak{n}-$\emph{cohomology groups of} $M$ are the $\mathfrak{l}-$%
modules obtained as the corresponding derived functors. These $\mathfrak{l}-$%
modules can be calculated by applying the functor 
\begin{equation*}
\text{Hom}_{U(\mathfrak{n})}(-,M)
\end{equation*}%
to the standard resolution of $\mathbb{C}$, this time by free left $U(%
\mathfrak{n})-$modules. In a natural way, one obtains a complex of $%
\mathfrak{l}-$modules and pth cohomology of this complex realizes the pth $%
\mathfrak{n}-$cohomology group

\begin{equation*}
H^{\text{p}}(\mathfrak{n},M).
\end{equation*}

\noindent Let $\mathfrak{n}^{\ast }$denote the $\mathfrak{l}-$module dual to 
$\mathfrak{n}$. Then, using the standard complexes and the natural
isomorphism of $\mathfrak{l}-$modules 
\begin{equation*}
\Lambda ^{p}\mathfrak{n}^{\ast }\otimes M\cong \text{Hom}(\Lambda ^{p}%
\mathfrak{n},M)
\end{equation*}%
one can deduce the following well known fact \cite[Section 2]{HS}:

\begin{proposition}
Suppose $M$ is a $\mathfrak{g-}$module. Let $\mathfrak{p}\subseteq \mathfrak{%
g}$ be a parabolic subalgebra with nilradical $\mathfrak{n}$ and Levi factor 
$\mathfrak{l}$. \newline
(a) Let $M^{\ast }$ denote the $\mathfrak{g}-$module dual to $M.$ Then there
is a natural isomorphism 
\begin{equation*}
H_{\text{p}}(\mathfrak{n},M^{\ast })\cong H^{\text{p}}(\mathfrak{n},M)^{\ast
}
\end{equation*}%
where $H^{\text{p}}(\mathfrak{n},M)^{\ast }$ denotes the $\mathfrak{l-}$%
module dual to $H^{\text{p}}(\mathfrak{n},M)$.\newline
(b) Let d denote the dimension of $\mathfrak{n}$. Then there is a natural
isomorphism 
\begin{equation*}
H_{\text{p}}(\mathfrak{n},M)\cong H^{\text{d-p}}(\mathfrak{n},M)\otimes
\Lambda ^{\text{d}}\mathfrak{n}\newline
\end{equation*}
\end{proposition}

\section{Harish-Chandra modules and globalizations}

\noindent Fix a maximal compact subgroup $K_{0}$ of $G_{0}$. Suppose we have
a linear action of $K_{0}$ on a complex vector space $M$. A vector $m\in M$
is called $K_{0}-$\emph{finite} if the span of the $K_{0}-$orbit of $m$ is
finite-dimensional and if the action of $K_{0}$ in this subspace is
continuous. The linear action of $K_{0}$ in $M$ is called $K_{0}-$\emph{%
finite }when every vector is $K_{0}-$finite. By definition, \emph{a
Harish-Chandra module} for $G_{0}$ is a finite length $\mathfrak{g}-$module $%
M$ equipped with a compatible, $K_{0}-$finite, linear action. One knows that
an irreducible $K_{0}-$module has finite multiplicity in a Harish-Chandra
module. For our purposes, it will also be useful to refer to a category of 
\emph{good }$K_{0}-$modules. A \emph{good }$K_{0}-$module will mean a
locally finite module such that each irreducible $K_{0}-$module has finite
multiplicity therein.

\smallskip

\noindent A representation of $G_{0}$ in a complete locally convex
topological vector space $V$ is called \emph{admissible} if $V$ has finite
length (with respect to closed invariant subspaces) and if each irreducible $%
K_{0}-$module has finite multiplicity in $V$. When $V$ is admissible, then
each $K_{0}-$finite vector in $V$ is differentiable and the subspace of $%
K_{0}-$finite vectors is a Harish-Chandra module. The representation is
called \emph{smooth} if every vector in $V$ is differentiable. In this case, 
$V$ is a $\mathfrak{g}-$module. For example, one knows that an admissible
representation in a Banach space is smooth if and only if the representation
is finite-dimensional.

\smallskip \smallskip

\noindent Given a Harish-Chandra module $M$, \emph{a globalization} $M_{%
\text{glob}}$ \emph{of }$M$ is an admissible representation of $G_{0}$ whose
underlying $(\mathfrak{g},K_{0})-$module of $K_{0}-$finite vectors is
isomorphic to $M$. By now, four canonical globalizations of Harish-Chandra
modules are known to exist. These are: the smooth globalization of Casselman
and Wallach \cite{C}, its dual (called: the distribution globalization),
Schmid's minimal globalization \cite{S} and its dual (the maximal
globalization). All four globalizations are smooth and functorial. In this
article we focus on the minimal and maximal globalizations of Schmid.

\smallskip \smallskip

\noindent The \emph{minimal globalization }$M_{\text{min}}$ of a
Harish-Chandra module $M$ is uniquely characterized by the property that any 
$(\mathfrak{g},K_{0})-$equivariant linear of $M$ onto the $K_{0}-$finite
vectors of an admissible representation $V$ lifts to a unique, continuous $%
G_{0}-$equivariant linear map of $M_{\text{min}}$ into $V$. In particular, $%
M_{\text{min}}$ embeds $G_{0}-$equivariantly and continuously into any
globalization of $M$. The construction of the minimal globalization shows
that it's realized on a \emph{DNF space}. This means that its continuous
dual, in the strong topology, is a nuclear Frech\'{e}t space. One knows that 
$M_{\text{min}}$ consists of analytic vectors and that it surjects onto the
analytic vectors in a Banach space globalization. Like each of the canonical
globalizations, the minimal globalization is functorially exact. In
particular, a closed $G_{0}-$invariant subspace of a minimal globalization
is the minimal globalization of its underlying Harish-Chandra module and a
continuous $G_{0}-$equivariant linear map between minimal globalizations has
closed range.

\smallskip

\noindent To characterize the maximal globalization, we introduce the $%
K_{0}- $finite dual on the category of Harish-Chandra modules. In
particular, let $M $ be a Harish-Chandra module. Then the algebraic dual $%
M^{\ast }$ of $M$ is a $\mathfrak{g}-$module and a $K_{0}-$module, but in
general not $K_{0}-$finite. We define $M^{\vee }$, \emph{the} $K_{0}-$\emph{%
finite (or Harish-Chandra) dual to} $M$, to be the subspace of $K_{0}-$%
finite vectors in $M^{\ast }$. Thus $M^{\vee }$ is a Harish-Chandra module.
In fact, the functor $M\mapsto M^{\vee }$ is exact on the category of good $%
K_{0}-$modules. We also have the formula 
\begin{equation*}
\left( M^{\vee }\right) ^{\vee }\cong M\text{.}
\end{equation*}%
The maximal globalization $M_{\text{max}}$ of $M$ can be defined by the
equation%
\begin{equation*}
M_{\text{max}}=\left( \left( M^{\vee }\right) _{\text{min}}\right) ^{\prime }
\end{equation*}%
where the last prime denotes the continuous dual equipped with the strong
topology. In particular, $M_{\text{max}}$ is a globalization of $M$. Observe
that the maximal globalization is an exact functor, since all functors used
in the definition are exact. Because of the minimal property of $M_{\text{min%
}}$, it follows that any globalization of $M$ embeds continuously and
equivariantly into $M_{\text{max}}$. Note that the continuous dual of a
maximal globalization is the minimal globalization of the dual
Harish-Chandra module.

\section{A conjecture by Vogan}

\noindent In order to introduce Vogan's conjecture, we need to be more
specific about the parabolic subalgebras we consider. Suppose $\mathfrak{p}$
is a parabolic subalgebra of $\mathfrak{g}$. We say that $\mathfrak{p}$ is 
\emph{nice} if $\mathfrak{g}_{0}\cap \mathfrak{p}=\mathfrak{l}_{0}$ is the
real form of a Levi factor $\mathfrak{l}$ of $\mathfrak{p}$. In this case $%
\mathfrak{l}$ is called \emph{the stable Levi factor}. When $\mathfrak{p}$
is nice, then every $G_{0}-$conjugate of $\mathfrak{p}$ is nice.

\smallskip

\noindent Suppose $\mathfrak{p}$ is nice and $\mathfrak{l}$ is the stable
Levi factor. Then we define \emph{the associated Levi subgroup} $L_{0}$ to
be the normalizer of $\mathfrak{p}$ in $G_{0}$. One knows that $L_{0}$ is a
real reductive group of Harish-Chandra class with Lie algebra $\mathfrak{l}%
_{0}$. Let 
\begin{equation*}
\theta :G_{0}\rightarrow G_{0}
\end{equation*}%
be a Cartan involution with fixed point set $K_{0}$. The parabolic
subalgebra will be called \emph{very nice} if $\theta (L_{0})=L_{0}$. In
this case $K_{0}\cap L_{0}$ is a maximal compact subgroup of $L_{0}$. One
knows that a nice parabolic subalgebra is $G_{0}-$conjugate to a very nice
parabolic subalgebra and that two very nice parabolic subalgebras are
conjugate under $K_{0}$ if and only if they are conjugate under $G_{0}$.

\bigskip

\noindent Throughout the remainder of this discussion, when $\mathfrak{p}$
is a very nice parabolic subalgebra, then $\mathfrak{n}$ will denote the
nilradical of $\mathfrak{p}$, $\mathfrak{l}$ will denote the stable Levi
factor and $L_{0}$ will denote the associated Levi subgroup. We fix the
maximal subgroup $K_{0}\cap L_{0}$ of $L_{0}$ and speak of Harish-Chandra
modules for $L_{0}$ accordingly. We have the following result \cite[%
Proposition 2.24]{HS}.

\begin{proposition}
Suppose $\mathfrak{p}$ is a very nice parabolic subalgebra and let $M$ be a
Harish-Chandra module for $G_{0}$. Then the $\mathfrak{n}-$homology groups
and $\mathfrak{n}-$cohomology groups of $M$ are Harish-Chandra modules for $%
L_{0}$.
\end{proposition}

\noindent Vogan's conjecture is the following.

\begin{conjecture}
Suppose $\mathfrak{p}$ is a very nice parabolic subalgebra and let $M$ be a
Harish-Chandra module for $G_{0}$. Suppose $M_{\text{glob}}$ indicates one
of the four canonical globalizations of $M$. Then the induced topologies of
the $\mathfrak{n}-$cohomology groups $H^{\text{p}}(\mathfrak{n},M_{\text{glob%
}})$ are Hausdorff and there are natural isomorphisms of $L_{0}-$%
representations 
\begin{equation*}
H^{\text{p}}(\mathfrak{n},M_{\text{glob}})\cong H^{\text{p}}(\mathfrak{n}%
,M)_{\text{glob}}
\end{equation*}%
where $H^{\text{p}}(\mathfrak{n},M)_{\text{glob}}$ denotes the corresponding
canonical globalization to $L_{0}$ of the Harish-Chandra module $H^{\text{p}%
}(\mathfrak{n},M)$. .
\end{conjecture}

\noindent Proposition 2.1 (b) shows that Vogan's conjecture holds for the $%
\mathfrak{n}-$cohomology groups if and only if it holds for the $\mathfrak{n}%
-$homology groups.

\smallskip

\noindent The conjecture is known to be true for the case of the minimal
globalization \cite{Br1}. We will now observe that Vogan's conjecture for
the dual representation is in fact equivalent to a certain purely algebraic
statement about $\mathfrak{n}-$homology groups and the Harish-Chandra dual
of a Harish-Chandra module.

\begin{proposition}
Suppose that $\mathfrak{p}$ is a very nice parabolic subalgebra and suppose
that Vogan's conjecture holds for the $\mathfrak{n}-$cohomology groups of
one of the four canonical globalizations. In particular, when $M$ is a
Harish-Chandra module, let $M_{\text{glob}}$ denote the globalization for
which Vogan's conjecture holds and let $M^{\text{glob}}$ denote the dual
globalization. Thus 
\begin{equation*}
M^{\text{glob}}\cong \left( \left( M^{\vee }\right) _{\text{glob}}\right)
^{\prime }\text{.}
\end{equation*}%
Then, in a natural way, the $\mathfrak{n}-$homology group $H_{\text{p}}(%
\mathfrak{n},M^{\text{glob}})$ is isomorphic to the dual globalization of $%
H^{\text{p}}(\mathfrak{n},M^{\vee })^{\vee }$. That is: 
\begin{equation*}
H_{\text{p}}(\mathfrak{n},M^{\text{glob}})\cong \left( H^{\text{p}}(%
\mathfrak{n},M^{\vee })^{\vee }\right) ^{\text{glob}}
\end{equation*}%
In particular, Vogan's conjecture holds for the dual globalization if and
only if there are natural isomorphisms 
\begin{equation*}
H_{\text{p}}(\mathfrak{n},M^{\vee })\cong H^{\text{p}}(\mathfrak{n},M)^{\vee
}
\end{equation*}%
for each p.
\end{proposition}

\noindent \textbf{Proof: }We assume the conjecture holds for $M_{\text{glob}%
} $. Since the continuous dual is exact on the category obtained by applying
the canonical globalization to Harish-Chandra modules, it follows, as in the
proof of Proposition 2.1, that 
\begin{equation*}
H_{\text{p}}(\mathfrak{n},\left( M_{\text{glob}}\right) ^{\prime })\cong H^{%
\text{p}}(\mathfrak{n},M_{\text{glob}})^{\prime }.
\end{equation*}%
Since $M^{\text{glob}}$ is given by $(\left( M^{\vee }\right) _{\text{glob}%
})^{\prime }$ it follows that%
\begin{eqnarray*}
H_{\text{p}}(\mathfrak{n},M^{\text{glob}}) &\cong &H_{\text{p}}\left( 
\mathfrak{n},\left( \left( M^{\vee }\right) _{\text{glob}}\right) ^{\prime
}\right) \cong H^{\text{p}}\left( \mathfrak{n},\left( M^{\vee }\right) _{%
\text{glob}}\right) ^{\prime }\cong \left( (H^{\text{p}}\left( \mathfrak{n}%
,\left( M^{\vee }\right) \right) _{\text{glob}}\right) ^{\prime } \\
&\cong &\left( H^{\text{p}}\left( \mathfrak{n},M^{\vee }\right) ^{\vee
}\right) ^{\text{glob}}\text{ \ \ \ }\blacksquare
\end{eqnarray*}

\noindent In this article we will show there are natural isomorphisms 
\begin{equation*}
H_{\text{p}}(\mathfrak{n},M^{\vee })\cong H^{\text{p}}(\mathfrak{n},M)^{\vee
}
\end{equation*}%
for $\mathfrak{p}$ a very nice parabolic subalgebra. We call this
isomorphism \emph{the} \emph{algebraic version of Vogan's conjecture}.

\section{The Natural Map and the Hochschild-Serre Spectral Sequence}

\noindent Through out the remainder of the discussion we fix a very nice
parabolic subalgebra $\mathfrak{p}$. Suppose $M$ is a Harish-Chandra module
for $G_{0}$. Then then the natural inclusion 
\begin{equation*}
M^{\vee }\rightarrow M^{\ast }
\end{equation*}%
induces a map 
\begin{equation*}
H_{\text{p}}(\mathfrak{n},M^{\vee })\rightarrow H_{\text{p}}(\mathfrak{n}%
,M^{\ast })\cong H^{\text{p}}(\mathfrak{n},M)^{\ast }\text{.}
\end{equation*}%
Since $H_{\text{p}}(\mathfrak{n},M^{\vee })$ is a Harish-Chandra module for $%
L_{0}$, it follows that the image of this map lies inside $H^{\text{p}}(%
\mathfrak{n},M)^{\vee }$. Our point is to show the resulting natural map 
\begin{equation*}
H_{\text{p}}(\mathfrak{n},M^{\vee })\rightarrow H^{\text{p}}(\mathfrak{n}%
,M)^{\vee }
\end{equation*}%
is an isomorphism.

\medskip

\noindent Before giving a general argument, we first make the following
observation. Suppose that $K_{0}\subseteq L_{0}$ (this is what happens, for
example when $G_{0}=SL(2,\mathbb{R})$). Then the standard resolution is a
complex of good $K_{0}-$modules (Section 3). Since the $K_{0}-$finite dual
is an exact functor on the category of good $K_{0}-$modules, one can argue
directly, as in the case for the ordinary dual that 
\begin{equation*}
H_{\text{p}}(\mathfrak{n},M^{\vee })\cong H^{\text{p}}(\mathfrak{n},M)^{\vee
}\text{.}
\end{equation*}%
Our general proof builds on this observation by introducing the
Hochschild-Serre spectral sequence \cite{KV}. In particular, let $\mathfrak{t%
}\subseteq $ $\mathfrak{g}$ be the complexified Lie algebra of $K_{0}$ and
let $\mathfrak{s}\subseteq \mathfrak{g}$ be the complexification of the
negative one eigenspace space for the Cartan involution $\theta $. Thus 
\begin{equation*}
\mathfrak{g}=\mathfrak{k}+\mathfrak{s}
\end{equation*}%
complexifies the Cartan decomposition for $\mathfrak{g}_{0}$. Since $%
\mathfrak{p}$ is $\theta -$ stable 
\begin{equation*}
\mathfrak{n}=\mathfrak{n}\cap \mathfrak{k}\oplus \mathfrak{n}\cap \mathfrak{s%
}\text{.}
\end{equation*}%
Observe that $\mathfrak{n}\cap \mathfrak{k}$ acts on $\mathfrak{n}\cap 
\mathfrak{s}$ by the adjoint representation. Roughly speaking, the
Hochschild-Serre spectral sequence gives us canonical ways to relate $H_{%
\text{q}}(\mathfrak{n}\cap \mathfrak{k},M^{\vee })\otimes \wedge \left( 
\mathfrak{n}\cap \mathfrak{s}\right) $ to $H_{\text{p}}(\mathfrak{n},M^{\vee
})$ and $H^{\text{q}}(\mathfrak{n}\cap \mathfrak{k},M)\otimes \wedge \left( 
\mathfrak{n}\cap \mathfrak{s}\right) ^{\ast }$ to $H^{\text{p}}(\mathfrak{n}%
,M)$. We will use this to deduce the desired result.

\medskip

\noindent We follow the development in \cite[Chapter V, Section 10]{KV} and
begin by selecting in $\wedge \left( \mathfrak{n}\cap \mathfrak{s}\right) $
a sequence $\left( V_{\text{p}}\right) _{\text{p=0}}^{\text{N}}$ of $%
K_{0}\cap L_{0}-$invariant subspaces such that

(a) $\wedge \left( \mathfrak{n}\cap \mathfrak{s}\right) =\oplus _{\text{p=0}%
}^{\text{N}}V_{\text{p}}$

(b) $V_{0}=\wedge ^{0}\left( \mathfrak{n}\cap \mathfrak{s}\right) =\mathbb{C}
$ and $V_{\text{N}}=\wedge ^{\text{R}}\left( \mathfrak{n}\cap \mathfrak{s}%
\right) $ where R=$\dim \left( \mathfrak{n}\cap \mathfrak{s}\right) $.

(c) There is a monotone increasing function $r($p$)\leq $p such that $V_{%
\text{p}}\subseteq \wedge ^{r(\text{p})}\left( \mathfrak{n}\cap \mathfrak{s}%
\right) $

(d) $\left( \mathfrak{n}\cap \mathfrak{k}\right) \cdot V_{\text{p}}\subseteq
\oplus _{\text{k=0}}^{\text{p-1}}V_{\text{k}}$.

\noindent The spectral sequences we need are as follows and can be phrased
in terms of a $\mathfrak{g}-$module $M$.

(i) There is a convergent spectral sequence 
\begin{equation*}
E_{\text{r}}^{\text{p,q}}\implies H_{\text{p+q}}(\mathfrak{n},M)
\end{equation*}%
with $E_{1}$ term 
\begin{equation*}
E_{1}^{\text{p,q}}=H_{\text{p+q}-r(\text{p})}(\mathfrak{n}\cap \mathfrak{k}%
,M)\otimes V_{\text{p}}.
\end{equation*}%
The differential $d_{\text{r}}$ has bidegree $(-$r,r$-1)$ and is a $%
K_{0}\cap L_{0}$ map when $M$ is a $K_{0}\cap L_{0}-$module.

(ii) There is a convergent spectral sequence 
\begin{equation*}
E_{\text{r}}^{\text{p,q}}\implies H^{\text{p+q}}(\mathfrak{n},M)
\end{equation*}%
with $E_{1}$ term 
\begin{equation*}
E_{1}^{\text{p,q}}=H^{\text{p+q}-r(\text{p})}(\mathfrak{n}\cap \mathfrak{k}%
,M)\otimes \left( V_{\text{p}}\right) ^{\ast }.
\end{equation*}%
The differential $d_{\text{r}}$ has bidegree $($r,$1.-$r$)$ and is a $%
K_{0}\cap L_{0}$ map when $M$ is a $K_{0}\cap L_{0}-$module. .

\medskip

\noindent We are now ready to prove the algebraic version of Vogan's
conjecture.

\begin{theorem}
Suppose $\mathfrak{n}$ is the nilradical of a very nice parabolic subalgebra
and let $M$ be a Harish-Chandra module. Then the natural map%
\begin{equation*}
H_{\text{p}}(\mathfrak{n},M^{\vee })\rightarrow H^{\text{p}}(\mathfrak{n}%
,M)^{\vee }
\end{equation*}%
is an isomorphism.
\end{theorem}

\noindent \textbf{Proof: }The proof is reminiscent of the proof given for 
\cite[Corollary 5.141]{KV}. First consider the spectral sequence associated
to the object $H_{\text{p}}(\mathfrak{n},M^{\ast })\cong H^{\text{p}}(%
\mathfrak{n},M)^{\ast }$. In particular, by dualizing everything in sight in
the spectral sequence for $H^{\text{p}}(\mathfrak{n},M)$ we obtain a
spectral sequence naturally isomorphic to the spectral sequence for $H_{%
\text{p}}(\mathfrak{n},M^{\ast })$ and with $E_{1}$ term 
\begin{equation*}
E_{1}^{\text{p,q}}(H_{\bullet }(\mathfrak{n},M^{\ast }))=H^{\text{p+q}-r(%
\text{p})}(\mathfrak{n}\cap \mathfrak{k},M)^{\ast }\otimes V_{\text{p}}.
\end{equation*}%
Thus the space of $K_{0}\cap L_{0}-$finite vectors in this term is given by%
\begin{equation*}
H^{\text{p+q}-r(\text{p})}(\mathfrak{n}\cap \mathfrak{k},M)^{\vee }\otimes
V_{\text{p}}.
\end{equation*}%
Next we show that this object is naturally isomorphic to the $E_{1}$ term
for the spectral sequence associated to $H_{\text{p}}(\mathfrak{n},M^{\vee
}) $. Letting $\widehat{K_{0}}$ denote the unitary dual of the group $K_{0}$%
, we write 
\begin{equation*}
M=\oplus _{\pi \in \widehat{K_{0}}}m(\pi )V_{\pi }
\end{equation*}%
where $V_{\pi }$ is a copy of the irreducible representation corresponding
to $\pi \in \widehat{K_{0}}$ and $m(\pi )$ is the multiplicity of $\pi $ in $%
M$. Thus we have 
\begin{equation*}
H_{\text{p+q}-r(\text{p})}(\mathfrak{n}\cap \mathfrak{k},M^{\vee })=H_{\text{%
p+q}-r(\text{p})}\left( \mathfrak{n}\cap \mathfrak{k},\left( \oplus _{\pi
\in \widehat{K_{0}}}m(\pi )V_{\pi }\right) ^{\vee }\right) =
\end{equation*}%
\begin{equation*}
H_{\text{p+q}-r(\text{p})}\left( \mathfrak{n}\cap \mathfrak{k},\oplus _{\pi
\in \widehat{K_{0}}}m(\pi )V_{\pi }^{\ast }\right) =\oplus _{\pi \in 
\widehat{K_{0}}}m(\pi )H_{\text{p+q}-r(\text{p})}\left( \mathfrak{n}\cap 
\mathfrak{k},V_{\pi }^{\ast }\right) =
\end{equation*}%
\begin{equation*}
\oplus _{\pi \in \widehat{K_{0}}}m(\pi )H^{\text{p+q}-r(\text{p})}\left( 
\mathfrak{n}\cap \mathfrak{k},V_{\pi }\right) ^{\ast }=\left( \oplus _{\pi
\in \widehat{K_{0}}}m(\pi )H^{\text{p+q}-r(\text{p})}\left( \mathfrak{n}\cap 
\mathfrak{k},V_{\pi }\right) \right) ^{\vee }=
\end{equation*}%
\begin{equation*}
H^{\text{p+q}-r(\text{p})}\left( \mathfrak{n}\cap \mathfrak{k},\oplus _{\pi
\in \widehat{K_{0}}}m(\pi )V_{\pi }\right) ^{\vee }=H^{\text{p+q}-r(\text{p}%
)}\left( \mathfrak{n}\cap \mathfrak{k},M\right) ^{\vee }\text{.}
\end{equation*}%
We now show how this leads to the desired result. In particular, let $E_{r}^{%
\text{p,q}}(H_{\bullet }(\mathfrak{n},M^{\vee }))$, $E_{r}^{\text{p,q}%
}(H_{\bullet }(\mathfrak{n},M^{\ast }))$ and $E_{r}^{\text{p,q}}(H^{\bullet
}(\mathfrak{n},M))$denote the $E_{\text{r}}$ terms of the corresponding
spectral sequences. Using induction, we want to see that 
\begin{equation*}
E_{r}^{\text{p,q}}(H_{\bullet }(\mathfrak{n},M^{\vee }))\cong E_{r}^{\text{%
p,q}}(H_{\bullet }(\mathfrak{n},M^{\ast }))_{K_{0}\cap L_{0}}
\end{equation*}%
where $E_{r}^{\text{p,q}}(H_{\bullet }(\mathfrak{n},M^{\ast }))_{K_{0}\cap
L_{0}}$ indicates the corresponding space of $K_{0}\cap L_{0}-$ finite
vectors. Indeed, using the fact that the $E_{\text{r+1}}$ terms are given by
the homology of a complex 
\begin{equation*}
d_{\text{r}}:E_{\text{r}}\rightarrow E_{\text{r}}
\end{equation*}%
we can determine the result from the fact the complex associated to $%
H_{\bullet }(\mathfrak{n},M^{\ast })$ is obtained by dualizing the complex
associated to $H^{\bullet }(\mathfrak{n},M)$ and the fact that terms $E_{r}^{%
\text{p,q}}(H^{\bullet }(\mathfrak{n},M))$ are good $K_{0}\cap L_{0}-$%
modules. Specifically 
\begin{equation*}
E_{r}^{\text{p,q}}(H\left( _{\bullet }(\mathfrak{n},M^{\ast })\right)
_{_{K_{0}\cap L_{0}}}\cong E_{r}^{\text{p,q}}\left( H^{\bullet }(\mathfrak{n}%
,M)\right) _{_{K_{0}\cap L_{0}}}^{\ast }=E_{r}^{\text{p,q}}\left( H^{\bullet
}(\mathfrak{n},M)\right) ^{\vee }
\end{equation*}%
therefore the result follows by induction since the $K_{0}\cap L_{0}-$
finite dual is an exact functor on the category of good $K_{0}\cap L_{0}-$%
modules.

\medskip

\noindent Finally, to deduce the main result, \ we use the filtrations of $%
H_{\text{p}}(\mathfrak{n},M^{\vee })$ and $H_{\text{p}}(\mathfrak{n},M^{\ast
})$ given by the corresponding spectral sequences. Then we can conclude the
final result from the following analog to \cite[Lemma 5.142]{KV}.

\begin{lemma}
Let $A$ be a good $(\mathfrak{l},K_{0}\cap L_{0})-$module and suppose $B$ is
an $(\mathfrak{l},K_{0}\cap L_{0})-$module (we do not assume $B$ is
necessarily $K_{0}\cap L_{0}-$finite). Suppose each module is endowed with $(%
\mathfrak{l},K_{0}\cap L_{0})-$filtrations 
\begin{equation*}
\begin{array}{c}
A^{N}\supseteq A^{N-1}\supseteq \cdots \supseteq A^{0}\supseteq A^{-1}=0 \\ 
B^{N}\supseteq B^{N-1}\supseteq \cdots \supseteq B^{0}\supseteq B^{-1}=0%
\end{array}%
\end{equation*}%
and let $j:A\rightarrow B$ be an $(\mathfrak{l},K_{0}\cap L_{0})-$map such
that $j(A^{\text{p }})\subseteq B^{\text{p}}$ for each p. Suppose the
induced maps 
\begin{equation*}
j^{\text{p}}:A^{\text{p}}/A^{\text{p-1}}\rightarrow \left( B^{\text{p}}/B^{%
\text{p-1}}\right) _{K_{0}\cap L_{0}}
\end{equation*}%
are isomorphisms for each p. Then the map 
\begin{equation*}
j:A\rightarrow B_{k_{0}\cap L_{0}}
\end{equation*}%
is an isomorphism.
\end{lemma}

\noindent \textbf{Proof of Lemma: }Since the functor taking the space of $%
K_{0}\cap L_{0}-$finite vectors is left exact on an appropriately defined
category of $(\mathfrak{l},K_{0}\cap L_{0})-$modules, the lemma follows as
in the proof of \cite[Lemma 5.142]{KV} after applying the functor of $%
K_{0}\cap L_{0}-$finite vectors \ $\blacksquare $

\noindent\ 

\noindent We can therefore conclude our main result.

\begin{theorem}
Let $G_{0}$ be a reductive Lie group of Harish-Chandra class, $%
K_{0}\subseteq G_{0}$ a maximal compact subgroup and $\mathfrak{g}$ the
complexified Lie algebra of $G_{0}$. Suppose $\mathfrak{n}$ is the
nilradical of a very nice parabolic subalgebra $\mathfrak{p}$ of $\mathfrak{g%
}$. Let $L_{0}\subseteq G_{0}$ denote the associated Levi subgroup and let $%
M_{\text{max}}$ denote the maximal globalization of a Harish-Chandra module $%
M$. Then, in a natural way, the $\mathfrak{n}-$cohomology groups $H_{\text{p}%
}(\mathfrak{n},M_{\text{max}})$ are representations of $L_{0}$ and for each
p, there are canonical isomorphisms 
\begin{equation*}
H^{\text{p}}(\mathfrak{n},M_{\text{max}})\cong H^{\text{p}}(\mathfrak{n},M)_{%
\text{max}}.
\end{equation*}
\end{theorem}

\end{document}